\documentclass[10pt]{amsart}

\usepackage[T1]{fontenc}
\usepackage[latin1]{inputenc}
\usepackage{amsmath}
\usepackage{amssymb}
\usepackage{amsthm}
\usepackage{dsfont}
\usepackage{color}
\usepackage{graphicx}

\theoremstyle{plain}\newtheorem{theo}{Theorem}
\theoremstyle{plain}\newtheorem{cor}{Corollary}
\theoremstyle{definition}\newtheorem{rem}{Remark}
\theoremstyle{plain}\newtheorem{defi}{Definition}[section]
\theoremstyle{plain}\newtheorem{lem}[defi]{Lemma}
\theoremstyle{plain}\newtheorem{prop}[defi]{Proposition}
\theoremstyle{definition}
\newtheorem*{prf}{Proof}

\newcommand{\N}{{\mathbb{N}}}

\begin{document}

\title[The Convergence rate of the Gibbs sampler for the 2-D Ising model  via a geometric bound]{The Convergence rate of the Gibbs sampler for the 2-D Ising model  via a geometric bound}
\author[ Brice Franke\,\,\, and\,\,\, Amine Helali  ]{Brice Franke\,\,\, and\,\,\, Amine Helali  } 
\address{Laboratoire d'Algèbre G\'{e}om\'{e}trie et Th\'{e}orie Spectrale LR/11/ES-53, D\'{e}partement de  Math\'{e}matique, Facult\'e des Sciences de Sfax, Universit\'{e} de Sfax, 3000 Sfax, Tunisia}

\address{Laboratiore de Math\'{e}matiques de Bretagne Atlantique UMR 6205, UFR Sciences et Techniques,
Universit\'{e} de Bretagne Occidentale, 6 Avenue Le Gorgeu, CS 93837, 29238 Brest, cedex 3, France }
\email{ brice.franke@univ-brest.fr}
\email{ amine.helali@univ-brest.fr}


\begin{abstract}
We study the geometric bound introduced by Diaconis and Stroock $(1991)$ of the Gibbs sampler for the two-dimensional Ising model with free boundary condition. The obtained result generalizes the method proposed by  Shiu and Chen $(2015)$  from dimension one to dimension two. Furthermore we observe that the new bound  improves the result given by Ingrassia  $(1994)$.
\end{abstract}                        
\maketitle
\section*{Introduction}
The Ising model is the most basic model in statistical mechanics having non trivial interaction. It has  many  applications in pattern analysis, molecular biology and image analysis. The distribution of the two-dimensional square lattice Ising model is:
 $$ \pi(x) = \frac{1}{Z_T} \exp \Bigg\{  \tfrac{1}{T} \Big{(}{  \displaystyle\sum_{j=1}^{n} \displaystyle\sum_{i=1}^{n-1}  x^j_i x^j_{i+1}    +\displaystyle\sum_{i=1}^{n} \displaystyle\sum_{j=1}^{n-1}  x^j_i x^{j+1}_{i}   } \Big{)} \Bigg\}  \,\,\,\,\,\,\,\,\, \forall \,x = (x^j_i)_{1\leq i,j \leq n }  \in  \chi\,\,, $$
where $\chi=\{-1, 1\}^{n^2}$ is the state space, $T$ is  a positive real representing  the temperature and  $$Z_T =  \displaystyle\sum_{x  \in  \chi} \exp \Bigg\{ \tfrac{1}{T}\Big{(}{  \displaystyle\sum_{j=1}^{n} \displaystyle\sum_{i=1}^{n-1}  x^j_i x^j_{i+1}    +\displaystyle\sum_{i=1}^{n} \displaystyle\sum_{j=1}^{n-1}  x^j_i x^{j+1}_{i}   }\Big{)} \Bigg\}$$ is the normalizing constant. Monte Carlo Markov chain $MCMC$ method is a very useful technique to draw samples from the Ising model. The Gibbs sampler introduced by Geman and Geman  (see \cite{Geman and Geman}) and the Metropolis-Hastings algorithm introduced by  Metropolis et al. (see \cite{Metropolis et al}) and Hastings (see \cite{Hastings})  are the most popular Monte Carlo Markov chain methods used in this context. Those two algorithms use an aperiodic and irreducible Markov chain which is reversible with respect to the measure $\pi$.
The reversibility is contained in the following $detailed$ $balance$ $equation$
$$Q(x,y) = \pi(x) P(x,y) = \pi(y) P(y,x)  = Q(y,x) \,\, \,\,\,\,\,\,\,\, \forall   \,x,y\in \chi .$$ 
Under the previous conditions (aperiodicity, irreducibility and reversibility) the measure $\pi$ is the unique invariant measure for the matrix $P$. The eigenvalues of $P$ can be arranged as follows:
$$ 1= \beta_0 > \beta_1 \geq \beta_2 \geq \cdots \geq \beta_{| \chi | -1} > -1.$$
There are two known criteria to measure the convergence rate of a Markov chain: On one hand, many researcher as Sinclair   (see \cite{Sinclair}), Frigessi et al.   (see \cite{Frigessi et al1}),  Ingrassia (see \cite{Ingrassia}),  Chen et al.   (see \cite{Shiu et al}) and Chen and Hwang   (see \cite{Chen et Hwang})  use the asymptotic variance to study the convergence rate of the MCMC algorithms. On the other hand, Diaconis and Stroock  (see \cite{Diacons Stroock}) use the total variation distance to quantify the convergence of the MCMC algorithms to their stationary distribution. We recall this result in the following theorem:
\begin{theo} [\mbox{Diaconis and Stroock 1991}]
If $P$ is a reversible Markov chain with unique invariant measure $\pi$ and $P$ is irreducible then for all $x \in \chi$ and $k  \in \N$:
$$ 4 \|P^k(x, .)  - \pi\|_{var}^2 = \left(\displaystyle\sum_{y \in \chi }|P(x, y) - \pi(y)| \right)^2 \leq \frac{1-\pi(x)}{\pi(x)} (\beta^*)^{2k} $$
where $\beta^* = \max\{\beta_1, | \beta_{| \chi | -1} | \}$.
\end{theo}
Many authors as Sinclair and Jerrum (see \cite{Sinclair and Jerrum}), Diaconis and Stroock (see \cite{Diacons Stroock}) and Sinclair  (see \cite{Sinclair}) introduced bounds for the eigenvalues $\beta_1$ and $\beta_{|\chi|-1}$. Ingrassia  (see \cite{Ingrassia}) shows that the result of Diaconis and Stroock (see  \cite{Diacons Stroock}) leads to the  tightest bound. This motivates us to have a closer look at explicit bounds that one can obtain from the  approach of  Diaconis and Stroock.\\
In order to do so, let us first remind the classical result of Diaconis and Stroock  (see \cite{Diacons Stroock}):\\
Let $\mathcal{G}(P)=(\chi, \, E)$ be the graph constructed with the dynamic $P$ where the state space $\chi$ is the vertex set and $E=\{(x, y)\mid P(x, y) >0 \}$ is the set of edges. Then, for each pair of distinct configurations $x$ and $y$ we choose a path $\gamma_{xy}$ in $\mathcal{G}(P)$ linking $x$ to $y$ . The irreducibility  of the matrix $P$ guarantees  that such paths exist. Finally, we define the set $\Gamma=\{\gamma_{xy}:\,\, x, y \in \chi \}$. Note that only one path $\gamma_{xy}$ for each pair of configurations $x$, $y$ is chosen. The second largest eigenvalue $\beta_1$ is then bounded from above as follows:
 \begin{eqnarray}\label{Geometric bound}
\beta_1 &\leq& 1 - \frac{1}{\kappa}
\end{eqnarray}
where 
\begin{eqnarray}
\kappa = \underset{e\in E}\max \,\,\, Q(e)^{-1} \sum_{\gamma_{xy} \ni e} | \gamma_{xy} | \pi(x) \pi(y), \label{2}
\end{eqnarray}
and where  $| \gamma_{xy} | $ designates the length of the path $\gamma_{xy}$.\\
In equation (\ref{2}) the maximum is over all directed edges in the graph $\mathcal{G}(P)$ and the sum is over all paths from the set $\Gamma$ passing through the fixed edge $e$. It is clear that $\kappa$ measures the bottlenecks (charged edges) in the graph $\mathcal{G}(P)$ . We notice that a  small $\kappa$ gives a better result. So, on one hand we should choose the shortest paths to link some pair of configurations $x$, $y$ and on the other hand we must avoid that many paths pass through the same edge in order to obtain the tightest bound for $\beta_1$.\\
In this paper we  study  the Gibbs sampler which chooses a random coordinate to be updated according to the conditional probabilities given the other coordinates. \\
The associated matrix  for the two-dimensional square lattice Ising model is
$$
P(x, y)=  \left\{
    \begin{array}{ll}
     \frac{ 1 }{ n^2 }  \pi(y^j_i | x) \,\,\,\,\,\,\,\, \,\,\,\,\, \,\,\,\,\,\,\,\,\,\,\,\,\,\,\,\,\,\,\,\, \,\,\mbox{if} \,\,\,\,\,\mbox{ d}(x, y)=1 \,\,\,\, \mbox{and} \,\,\,\, x^j_i \neq y^j_i
\\    
\\
1- \frac{ 1 }{ n^2 }  \displaystyle\sum_{i,j=1}^{n}  \pi(y^j_i | x) \,\,\,\,\,\,\,\, \,\, \mbox{  if }\,\,\,\,\, \mbox{d}(x, y)=0
\\
\\
0 \,\,\,\,\,\,\,\,\,\,\,\,\,\,\,\,\,\,\,\,\,\,\,\,\,\,\,\,\,\,\,\,\,\,\,\,\,\,\,\,\,\,\,\,\,\,\,\,\,\,\,\,\,\,\,\, \mbox{ else}   
    \end{array}
\right.
$$
where 
\begin{itemize} 
\item $  \pi(y^j_i | x) =  \frac{  \pi(x^1_1,  \cdots , x^j_{i-1}, y^j_i,  x^j_{i+1}, \cdots, x^n_n  )}{  \pi(x^1_1,  \cdots , x^j_{i-1}, y^j_i,  x^j_{i+1}, \cdots, x^n_n  ) +  \pi(x^1_1,  \cdots , x^j_{i-1}, -y^j_i,  x^j_{i+1}, \cdots, x^n_n  )};$\\
\item $\mbox{d}(x, y)=\sharp \Big\{i,\,\,\, x^j_i \neq y^j_i \,\, \mbox{for} \,\, i,j= 1, \cdots, n  \Big\}$ designates the number of sites that differ between two configurations $x$ and $y$.
\end{itemize}
We apply the bound for $\beta_1$  introduced by Diaconis and Stroock  (see \cite{Diacons Stroock}) to the Gibbs sampler for the   two-dimensional Ising model with two states. The computation is based on some method introduced by Shiu and Chen  (see \cite{Shiu chen}). In their paper they treat the one-dimensional lattice case. \\
The one-dimensional lattice case with multiple states which is also called Potts model was investigated by Helali with the same techniques (see \cite{Helali}).\\
In the first section, we first define a path for each pair of configurations $(x, y)$ from the state space $ \chi=\{-1, +1\}^{n^2}$. Then, we turn to compute explicitly the expression:  
\begin{align}
\displaystyle \max_{e \in E}Q(e)^{-1} \sum_{(x, y):\,\,\gamma_{xy} \ni e} | \gamma_{xy} | \pi(x) \pi(y).\label{bounda}
\end{align}
In Proposition $1.1$, we give an upper bound for the expression in (\ref{bounda}) for different classes  of edges (interior, corners, etc $\cdots$) in the square lattice $\{1, \cdots, n\}^2$. We notice that it is difficult to complete exact computations for the bounds given in Proposition $1.1$. In Theorem $2$ we present a bound for $\beta_1$ which results from some rough estimation of the terms from Proposition $1.1$. To be able to use Theorem $1$ we are referring to Ingrassia (see Theorem $5.3$ in \cite{Ingrassia}) which gives a lower bound for the smallest eigenvalue $\beta_{|\chi|-1}$. The main theorem is obtained by the fact that $|\beta_{|\chi|-1}|$ is smaller than the upper bound of  $\beta_1$. 
In the third  section we compare our result with existing bounds from the literature. We notice that the main result of this paper generalizes the one introduced by Shiu and Chen  in \cite{Shiu chen} to higher dimension. It also improves the result of Ingrassia   (see \cite{Ingrassia}). The last section contains the proofs for the main results
 \section{Main result}
\subsection{Selection of paths}
 To be able to use the result of Diaconis and Stroock (see \cite{Diacons Stroock}) in the computation of the bound for the second largest eigenvalue, we should fix a collection of paths linking any $x$ to any $y$ from $\chi$. To get the best possible result, we should use the shortest paths $\gamma_{xy}$ to obtain smallest possible  $\kappa$ which  gives smaller upper bound for $\beta_1$. 
An edge of the graph $\mathcal{G}(P)$ takes the following form: $e=(e^-, e^+)$  where
 $$e^-=\begin{pmatrix}
   z^1_1 & \cdots   & \cdots  & z^1_p& \cdots &\cdots&z^1_n  \\
      \vdots  &     &  &\vdots &  & &\vdots   \\ 
   z^{q-1}_1 & \cdots  & z^{q-1}_{p-1}& z^{q-1}_p & z^{q-1}_{p+1} &\cdots&   z^{q-1}_n \\    
    z^q_1 &\cdots  & z^q_{p-1}& z^q_p & z^q_{p+1} &\cdots&   z^{q}_n\\   
      z^{q+1}_1 & \cdots  & z^{q+1}_{p-1}& z^{q+1}_p & z^{q+1}_{p+1} &\cdots&   z^{q+1}_n \\   
  \vdots    &   & & \vdots &   && \vdots   \\ 
   z^n_1& \cdots &  \cdots & z^n_p & \cdots  &\cdots& z^n_n
\end{pmatrix} 
$$
and  $e^+$ takes the same form except that we have $-z_p^q$ in the $(p, q)-th$ position.\\
Without loss of generality we rearrange the lattice as a vector $x=(x_1, \cdots, x_{n^2})$ then we use the same kind of paths introduced by Shiu and Chen in \cite{Shiu chen}. For a given pair $(x, y) \in \chi^2$ there exist an increasing sequence  $d_1, \cdots, d_m \in \{1, \cdots, n^2  \} $ such that:  
\begin{itemize}
\item $ x_i \neq y_i $ for $i \in \{ d_1, \cdots, d_m \}$  
\vspace{0.3cm}
\item $x_i = y_i $  otherwise.
\end{itemize}
The path linking $x$ to $y$ is defined as follows:
\begin{align}
(x_1, \cdots, x_{n^2}) &= (y_1, \cdots, y_{d_1-1}, x_{d_1}, x_{d_1+1}, \cdots, x_{d_2-1}, x_{d_2}, x_{d_2+1}, \cdots, x_{n^2}) \nonumber \\
&\rightarrow  (y_1, \cdots, y_{d_1-1}, y_{d_1}, x_{d_1+1}, \cdots, x_{d_2-1}, x_{d_2}, x_{d_2+1}, \cdots, x_{n^2}) \nonumber \\
&= (y_1, \cdots, y_{d_1-1}, x_{d_1}, y_{d_1+1}, \cdots, y_{d_2-1}, x_{d_2}, x_{d_2+1}, \cdots, x_{n^2}) \nonumber \\
&\rightarrow  (y_1, \cdots, y_{d_1-1}, y_{d_1}, y_{d_1+1}, \cdots, y_{d_2-1}, y_{d_2}, x_{d_2+1}, \cdots, x_{n^2}) \nonumber \\
&\vdots  \nonumber \\
& \rightarrow (y_1, \cdots, y_{n^2}).    \nonumber
\end{align}
Those paths subsequently update the differing sites in the configurations $x$ and $y$.
\subsection{Geometric bound for the second largest eigenvalue}     
In what follows we will find an upper bound for $$ \kappa = \underset{e \in E}\max \,\, Q(e)^{-1} \sum_{\gamma_{xy} \ni e} | \gamma_{xy} |  \pi(x) \pi(y).$$ 
 First, let $e=(e^-, e^+)$ be a fixed edge.  The flux in equilibrium  associated to this edge is:
$$Q(e)  = \pi(e^-)    P( e^-, e^+)$$     
and the transition matrix $P(e^-, e^+)$ is described in the following lemma:
\begin{lem}
Let $x$ be a configuration. From the position of the site who is to be updated we distinguish three principale cases:\\
1) If $(p, q)\in \{2, \cdots, n-1\}^2$, each configuration has four  neighbors and
 $$P( e^-, e^+) =  \frac{1}{n^2(1 + e^{\frac{2}{T}(z^q_{p-1}z^q_p + z^q_p z^q_{p+1} +z^{q-1}_{p}z^q_p + z^q_p z^{q+1}_{p} )})}.$$
2) If $(p, q) \in \{2, \cdots, n-1\} \times \{1\}$, each configuration has three  neighbors and
$$P( e^-, e^+)  =\frac{1}{n^2(1 + e^{\frac{2}{T}(  z^{q-1}_{1}z^{q}_{1} + z^{q}_{1}z^{q+1}_{1} + z^{q}_{1}z^{q}_{2})})}$$
and similar results hold for: \\
i) $(p, q)\in \{2, \cdots, n-1\} \times \{n\}$:
$$P( e^-, e^+)  =\frac{1}{n^2(1 + e^{\frac{2}{T}(  z^{q}_{n-1}z^{q}_{n} + z^{q-1}_{n}z^{q}_{n} + z^{q}_{n}z^{q+1}_{n})})},$$
ii) $(p, q)\in  \{1\} \times \{2, \cdots, n-1\}$:
$$P( e^-, e^+)  =\frac{1}{n^2(1 + e^{\frac{2}{T}(  z^{1}_{p-1}z^{1}_{p} + z^{1}_{p}z^{1}_{p+1} + z^{1}_{p}z^{2}_{p})})},$$
iii)  $(p, q)\in  \{n\} \times \{2, \cdots, n-1\}$:
$$P( e^-, e^+)  =\frac{1}{n^2(1 + e^{\frac{2}{T}(  z^{n}_{p-1}z^{n}_{p} + z^{n}_{p}z^{n}_{p+1} + z^{n}_{p}z^{n-1}_{p})})}.$$
3) If $(p, q)=(1, 1)$, each configuration has two  neighbors and
$$ P( e^-, e^+) = \frac{1}{n^2(1 + e^{\frac{2}{T}(z^1_{1}z^1_2 + z^1_1 z^2_{1})})}$$
and similar results hold for:\\ 
i) $(p, q)=(1, n)$:
$$ P( e^-, e^+) = \frac{1}{n^2(1 + e^{\frac{2}{T}(z^n_{1}z^n_2 + z^{n-1}_1 z^n_{1})})},$$
ii) $(p, q)=(n, 1)$:
$$ P( e^-, e^+) = \frac{1}{n^2(1 + e^{\frac{2}{T}(z^1_{n-1}z^1_{n} + z^1_n z^2_{n})})},$$
iii) $(p, q)=(n, n)$:
$$ P( e^-, e^+) = \frac{1}{n^2(1 + e^{\frac{2}{T}(z^n_{n-1}z^n_n + z^{n-1}_n z^n_{n})})}.$$
\end{lem}
The proof of this lemma is given in section $4$.\\
For each class of edges, we give an upper bound for $\kappa$ defined in equation (\ref{2})  in the following proposition:
\begin{prop}
The second largest eigenvalue of the Gibbs sampler for the two-dimensional  Ising model satisfies: $$ \beta_1 \leq 1 - \frac{1}{\kappa} \,\,\,\,\,\,\,\,  where \,\, \kappa \,\, is \,\, bounded \,\, as \,\, follow:   $$ 
1) For $(p, q)=(1, 1)$ or $(p, q)=(n, n)$  we have: 
\begin{align}
 Q(e)^{-1} \sum_{\gamma_{xy} \ni e} |\gamma_{xy}| \pi(x) \pi(y) &\leq \frac{n^4}{2}(1 + e^{\frac{4}{T}})  .\nonumber 
 \end{align}
2) For $(p, q)=(1, n)$    we have:
\begin{align}
 Q(e)^{-1} \sum_{\gamma_{xy} \ni e} |\gamma_{xy}| \pi(x) \pi(y) &\leq 2 n^4 e^{\frac{n-1}{T}} \sum_{w \in \chi:\,\, w_n^1=1} \pi(w)  \exp \bigg\{ \tfrac{1}{T} \Big{(} 2 (1 - w_{n}^{2}  )  \nonumber \\
 &+ \displaystyle \sum_{i=1}^{n-1} \big{(}  w_{i}^1  - w_{i}^{2}  -  w_{i}^{1} w_{i}^{2}    \big{)}    \Big{)} \bigg\}.\nonumber
\end{align}
and a similar result holds for  $(p, q)=(n, 1)$.\\
3) For $(p,q)\in \{2, \cdots, n-1\} \times \{1\} $  we have:  
\begin{align}
 Q(e)^{-1} \sum_{\gamma_{xy} \ni e} |\gamma_{xy}| \pi(x) \pi(y)&= n^4 e^{\frac{n+1}{T}} \sum_{w\in \chi: \,w_p^1=1 }  \pi(w)   \exp\bigg\{ \tfrac{1}{T} \Big{(} - 2 (w^{q}_{2} +w_1^{q+1})   + \displaystyle\sum_{i=2}^{n} \big{(}    w^{q-1}_i        \nonumber \\
& -  w^{q}_{i} -   w^{q-1}_i w^{q}_{i}  \big{)}  \Big{)}  \bigg\}    +  n^4 e^{\frac{n+1}{T}}  \sum_{w\in \chi: \,w_p^1=-1}  \pi(w)   \exp \bigg\{ \tfrac{1}{T} \Big{(} 2(  1  +    w^{q-1}_1  )\nonumber \\
&+  \displaystyle\sum_{i=2}^{n} \big{(}   w^{q-1}_i    - w^{q}_{i}  -   w^{q-1}_i w^{q}_{i}  \big{)}   \Big{)}   \bigg\}.\nonumber
\end{align}
 A similar result holds   for   $(p, q) \in \{2, \cdots, n-1 \}\times \{n\}$. Moreover with the same tricks we obtain an upper bound in the cases where
  $$(p, q)\in \{1\} \times \{ 2, \cdots, n-1\}\,\,\,\,\,\,  \mbox{or}\,\,\,\,\,\, (p, q)\in \{n\} \times \{ 2, \cdots, n-1\}.$$
4) For $(p, q) \in \{ 2, \cdots, n-1 \}^2$ we have:   
\begin{align}
 Q(e)^{-1} \sum_{\gamma_{xy} \ni e} |\gamma_{xy}| \pi(x) \pi(y) &\leq  2 n^4 e^{\frac{1}{T}(n-1)}  \sum_{w \in \chi: w^p_q=1} \pi(w)  \exp \bigg\{  \tfrac{1}{T}  \Big{(}  \displaystyle\sum_{i=1}^{p-1}  \big{(}   -w_i^q +   w^{q+1}_{i}  - w^q_i w^{q+1}_{i} \big{)}    \nonumber \\
&- 2(  w^q_{p+1} +      w^{q+1}_{p} ) +\displaystyle \sum_{i=p+1}^{n}  (    w^{q-1}_i  -  w^{q}_{i} -  w^{q-1}_i w^{q}_{i}    \big{)}    \Big{)} \bigg\}.  \nonumber 
\end{align}
\end{prop}
The proof of this proposition is given later in section $4$.\\
After this step, we notice that  it is difficult to give an exact value for the  sums  in the previous proposition. So, in what follows we give an upper bound for those sums in order  to obtain an upper bound for $\kappa$ defined in equation (\ref{2}). The main result is given in the following theorem:
\begin{theo}
 The second largest eigenvalue of the Gibbs sampler for the two-dimensional Ising model with two states satisfies:
 $$ \beta_1 \leq 1 - n^{-4} \displaystyle \exp{\big\{-\tfrac{2}{T}\left(2n+1\right)\big\}} . $$
\end{theo}
The proof of this theorem is given in section $4$.\\
To be able to use  Theorem $1$ given by Diaconis and Stroock \cite{Diacons Stroock} we must give an upper bound for the second largest eigenvalue in absolute value. Until now, we have given an upper bound for the second largest eigenvalue. So we turn now to control the smallest eigenvalue of the  Gibbs sampler for the two-dimensional Ising model.
\subsection{Bound of the second largest eigenvalue in absolute value}
Theorem  $5.3$ introduced by Ingrassia  in \cite{Ingrassia} gives a lower bound  for the smallest eigenvalue as:
$$ \beta_{| \chi | - 1} \geq -1 + \frac{2}{1+(c-1)e^{\frac{\Delta}{T}}}.$$
For the two-dimensional Ising model one has $c=2$ and $\Delta=4$. Therefore for large $n$:
$$| \beta_{| \chi | - 1} | \leq |  -1 + \frac{2}{1+ e^{\frac{4}{T}}} | = 1- \frac{2}{1+ e^{\frac{4}{T}}} < 1- \frac{2}{e^{\frac{4}{T}}+ e^{\frac{4}{T}}} < 1- e^{\frac{-4}{T}}< 1- n^{-4} e^{\frac{-2}{T}(2n+1)}.$$
An upper bound of $\beta^*$ is given in the following corollary: 

\begin{cor}
The second largest eigenvalue in absolute value of the Gibbs sampler for the two-dimensional Ising model with two states satisfies:
 $$ \beta^* \leq 1 - n^{-4} \displaystyle \exp{\big\{-\tfrac{2}{T}\left(2n+1\right)\big\}} . $$
\end{cor}
\begin{prf}
We give  an upper bound for $| \beta_{| \chi | - 1} |$ in the above computation. Then, combining this result with the result given in Theorem $2$  finishes the proof.
\end{prf}
\section{Comparison}
Ingrissia (see \cite{Ingrassia}) develop a method to give an upper bound for the second largest eigenvalue of the Gibbs sampler for the  general Ising model and he obtains: $$ \beta_1 \leq 1- \frac{Z_T}{ b_{\gamma}\, \Gamma_{\gamma} \, c \, |S| } e^{-\frac{m}{T} },$$
where $Z_T$ is the normalizing constant, $S$ is the lattice of sites,  $\Gamma$ the collection of paths, $\gamma_{\Gamma}$ the maximum length of each path $\gamma_{xy} \in \Gamma$,  $b_{\Gamma}$ is the maximum number of paths containing any edge of $\Gamma$, $c$ is the number of configurations that differ by only one site and $m$ is the least total elevation gain of the Hamiltonian function in the sense which is described by Holley and Stroock   (see \cite{Holley and Stroock}). \\
For the two-dimensional square lattice Ising model with two states we have:
 $$ \Gamma_{\gamma}=n^2,\,\, b_{\gamma}= 2^{n^2-1},\,\, c= 2,\,\, m=4,\,\, |S|=n^2\,\, \mbox{and}\,\, Z_T \leq 2 ( 1+ e^{-\frac{1}{2T}} )^{n^2-1}.$$ 
 Which leads to:
$$  \beta_1 \leq 1- n^{-4} e^{-\frac{4}{T}} \Big{(}\frac{1+e^{-\frac{1}{2T}}}{2}\Big{)}^{n^2-1}.$$  
The comparison of the two results amounts to compare $f(T)=e^{\frac{4}{T}}$  and  $g(T)=\frac{2}{1+e^{-\frac{1}{2T}}}$. The following figure represents the two graphs.
\begin{center}
\includegraphics[width= 11cm]{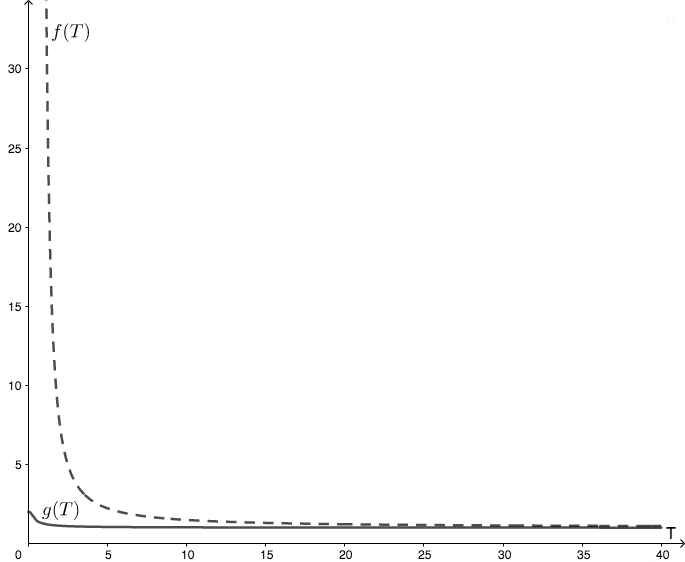} \\  
\vspace{0.2 cm}
\end{center}
We notice that for non-zero temperature, i.e.: $T > 0$ we have $e^{\frac{4}{T}}  \geq \frac{2}{1+e^{-\frac{1}{2T}}} $. However for all sufficiently large natural numbers $n$ holds: $$ \exp\Big\{\tfrac{2}{T}(2n-1) \Big\} \leq \exp\Big\{\tfrac{4n}{T} \Big\} = \left(e^{\frac{4}{T}} \right)^n \ll \left(\frac{2}{1+e^{-\frac{1}{2T}}}\right)^{n^2-1}.$$
This is quite natural to consider because MCMC method is used to give samples from probability measures defined on large spaces ($n \sim 10^{23}$). Then the main result of this paper improves the one introduced by Ingrassia (see \cite{Ingrassia}) for large $n$.
\section{ conclusion}
This paper deals with the bound of the second largest eigenvalue in absolute value of the Gibbs sampler for the  two-dimensional  Ising model with two states $\pm1$. The main result generalizes some methods of Shiu and Chen (see \cite{Shiu chen}) obtained in dimension one to dimension two. It also improves  Ingrassia's bound (see \cite{Ingrassia}).
\section{Proof of the main result}
\subsection{Proof of  lemma $1.1$}
We concentrate on the case where  $(p, q) \in \{ 2, \cdots, n-1 \}^2$. We have:
\begin{align}
\pi(e^-)&= \frac{1}{Z_T} \exp \bigg\{ \tfrac{1}{T} \Big{(} \displaystyle\sum_{j=1}^{q-1} \displaystyle\sum_{i=1}^{n-1}  x^j_i x^j_{i+1}  +  \displaystyle\sum_{i=1}^{p-2}  z^q_i z^q_{i+1} +  z^q_{p-1} z^q_{p} + z^q_{p} z^q_{p+1} +  \displaystyle\sum_{i=p+1}^{n-1}  z^q_i z^q_{i+1}      +\displaystyle\sum_{j=q+1}^{n} \displaystyle\sum_{i=1}^{n-1}  z^j_i z^j_{i+1}  \nonumber \\
&+  \displaystyle\sum_{i=1}^{p-1} \big{(} \displaystyle\sum_{j=1}^{q-2}  z^j_i z^{j+1}_{i} +    z^{q-1}_i z^{q}_{i} + z^q_i z^{q+1}_{i}   +   \displaystyle\sum_{j=q+1}^{n-1}  z^j_i z^{j+1}_{i}    \big{)}  +  \displaystyle\sum_{j=1}^{q-2}  z^j_p z^{j+1}_{p} +    z^{q-1}_p z^{q}_{p} + z^{q}_p z^{q+1}_{p}     \nonumber \\
&+   \displaystyle\sum_{j=q+1}^{n-1}  z^j_p z^{j+1}_{p}  +  \displaystyle\sum_{i=p+1}^{n} \big{(}  \displaystyle\sum_{j=1}^{q-2}  z^j_i z^{j+1}_{i}  +   z^{q-1}_i z^{q}_{i} +  z^q_i z^{q+1}_{i} +   \displaystyle\sum_{j=q+1}^{n-1}  z^j_i z^{j+1}_{i}    \big{)}  \Big{)}\bigg\} .   \nonumber
\end{align}
and $\pi(e^+)$ takes a similar form except in the $(p, q)-$th position which is equal to $-z_p^q$. This yields: 
\begin{align}
P(e^-, e^+)&= \frac{1}{n^2} \frac{\pi(e^+)}{\pi(e^+) + \pi(e^-)} \nonumber \\
&= \frac{ e^{\frac{-1}{T}(z^q_{p-1}z^q_p + z^q_p z^q_{p+1} +z^{q-1}_{p}z^q_p + z^q_p z^{q+1}_{p} )}  }{e^{\frac{1}{T}(z^q_{p-1}z^q_p + z^q_p z^q_{p+1} +z^{q-1}_{p}z^q_p + z^q_p z^{q+1}_{p} )}  + e^{\frac{-1}{T}(z^q_{p-1}z^q_p + z^q_p z^q_{p+1} +z^{q-1}_{p}z^q_p + z^q_p z^{q+1}_{p} )} }  \nonumber \\
&=   \frac{1}{ 1+ e^{\frac{2}{T}(z^q_{p-1}z^q_p + z^q_p z^q_{p+1} +z^{q-1}_{p}z^q_p + z^q_p z^{q+1}_{p} )} } .\nonumber 
\end{align}
With similar techniques we can compute the transition matrix for each class of edges and the results are presented in lemma $1.1$.
\subsection{Proof of  proposition $1.2$}
Depending on the position of the site who is to be updated we distinguish several cases.\\
Lets first assume that  $(p, q) \in \{ 2, \cdots, n-1 \}^2$. Then the transition matrix is:
$$P( e^-, e^+)= \frac{1}{n^2(1 + e^{\frac{2}{T}(z^q_{p-1}z^q_p + z^q_p z^q_{p+1} +z^{q-1}_{p}z^q_p + z^q_p z^{q+1}_{p} )})}.$$
According to the selection of the paths done in section $1.1.$, if a  pair $(x, y)$ satisfies $ \gamma_{xy} \ni e$, then the configurations $x$ and $y$ must have the following form:
$$
x=\begin{pmatrix}
   x^1_1 & \cdots   & \cdots  & x^1_p& \cdots &\cdots& x^1_n  \\
      \vdots&       &  &\vdots &  & & \vdots   \\ 
   x^{q-1}_1&  \cdots  & x^{q-1}_{p-1}& x^{q-1}_p & x^{q-1}_{p+1} &\cdots&   x^{q-1}_n \\    
    x^q_1 &\cdots  & x^q_{p-1}& z^q_p & z^q_{p+1} &\cdots &  z^{q}_n\\   
      z^{q+1}_1 & \cdots  & z^{q+1}_{p-1}& z^{q+1}_p & z^{q+1}_{p+1} &\cdots &  z^{q+1}_n \\   
  \vdots   & &    & \vdots &  & & \vdots   \\ 
   z^n_1& \cdots &  \cdots & z^n_p & \cdots  &\cdots& z^n_n
\end{pmatrix} $$
\,\,
\mbox{and} \,\,
$$y=\begin{pmatrix}
   z^1_1&  \cdots   & \cdots  & z^1_p& \cdots &\cdots& z^1_n  \\
      \vdots &      &  &\vdots &  & &\vdots   \\ 
   z^{q-1}_1 & \cdots  & z^{q-1}_{p-1}& z^{q-1}_p & z^{q-1}_{p+1} &\cdots  & z^{q-1}_n \\    
    z^q_1 &\cdots  & z^q_{p-1}& -z^q_p & y^q_{p+1} &\cdots &  y^{q}_n\\   
      y^{q+1}_1 & \cdots  & y^{q+1}_{p-1}& y^{q+1}_p & y^{q+1}_{p+1} &\cdots &  y^{q+1}_n \\   
  \vdots   &    & & \vdots &   && \vdots   \\ 
   y^n_1& \cdots &  \cdots & y^n_p & \cdots  &\cdots& y^n_n
\end{pmatrix} .
$$
We then have:
\begin{align}
\pi(x)&= \frac{1}{Z_T} \exp \bigg\{ \tfrac{1}{T} \Big{(} \displaystyle\sum_{j=1}^{q-1} \displaystyle\sum_{i=1}^{n-1}  x^j_i x^j_{i+1}  +  \displaystyle\sum_{i=1}^{p-2}  x^q_i x^q_{i+1} +  x^q_{p-1} z^q_{p} + z^q_{p} z^q_{p+1} +  \displaystyle\sum_{i=p+1}^{n-1}  z^q_i z^q_{i+1}      +\displaystyle\sum_{j=q+1}^{n} \displaystyle\sum_{i=1}^{n-1}  z^j_i z^j_{i+1}  \nonumber \\
&+  \displaystyle\sum_{i=1}^{p-1} \big{(} \displaystyle\sum_{j=1}^{q-2}  x^j_i x^{j+1}_{i} +    x^{q-1}_i x^{q}_{i} + x^q_i z^{q+1}_{i}   +   \displaystyle\sum_{j=q+1}^{n-1}  z^j_i z^{j+1}_{i}    \big{)}  +  \displaystyle\sum_{j=1}^{q-2}  x^j_p x^{j+1}_{p} +    x^{q-1}_p z^{q}_{p} + z^{q}_p z^{q+1}_{p}     \nonumber \\
&+   \displaystyle\sum_{j=q+1}^{n-1}  z^j_p z^{j+1}_{p}  +  \displaystyle\sum_{i=p+1}^{n} \big{(}  \displaystyle\sum_{j=1}^{q-2}  x^j_i x^{j+1}_{i}  +   x^{q-1}_i z^{q}_{i} +  z^q_i z^{q+1}_{i} +   \displaystyle\sum_{j=q+1}^{n-1}  z^j_i z^{j+1}_{i}    \big{)}  \Big{)}\bigg\}.   \nonumber
\end{align}
and
\begin{align}
 \pi(y)&= \frac{1}{Z_T} \exp \bigg\{ \tfrac{1}{T} \Big{(} \displaystyle\sum_{j=1}^{q-1} \displaystyle\sum_{i=1}^{n-1}  z^j_i z^j_{i+1}  +  \displaystyle\sum_{i=1}^{p-2}  z^q_i z^q_{i+1} -  z^q_{p-1} z^q_{p} - z^q_{p} y^q_{p+1} +  \displaystyle\sum_{i=p+1}^{n-1}  y^q_i y^q_{i+1}      +\displaystyle\sum_{j=q+1}^{n} \displaystyle\sum_{i=1}^{n-1}  y^j_i y^j_{i+1}  \nonumber \\
&+  \displaystyle\sum_{i=1}^{p-1} \big{(} \displaystyle\sum_{j=1}^{q-2}  z^j_i z^{j+1}_{i} +    z^{q-1}_i z^{q}_{i} + z^q_i y^{q+1}_{i}   +   \displaystyle\sum_{j=q+1}^{n-1}  y^j_i y^{j+1}_{i}    \big{)}  +  \displaystyle\sum_{j=1}^{q-2}  z^j_p z^{j+1}_{p} -    z^{q-1}_p z^{q}_{p} - z^{q}_p y^{q+1}_{p}     \nonumber \\
&+   \displaystyle\sum_{j=q+1}^{n-1}  y^j_p y^{j+1}_{p}  +  \displaystyle\sum_{i=p+1}^{n} \big{(}  \displaystyle\sum_{j=1}^{q-2}  z^j_i z^{j+1}_{i}  +   z^{q-1}_i y^{q}_{i} +  y^q_i y^{q+1}_{i} +   \displaystyle\sum_{j=q+1}^{n-1}  y^j_i y^{j+1}_{i}    \big{)}  \Big{)}\bigg\}.   \nonumber
\end{align}
Similarly, we can compute the expression of the measure $\pi(e^-)$. We then obtain:
\begin{align}
 Q(e)^{-1} \pi(x) \pi(y) &=  \frac{\pi(x) \pi(y)}{\pi(e^-)    P( e^-, e^+)}  \nonumber \\
             &=   \frac{n^2}{Z_T} \left( \exp \Big\{\tfrac{-2}{T} (z^q_{p-1} z^q_{p})\Big\}  + \exp \Big\{\tfrac{2}{T}(z^q_p z^q_{p+1} +z^{q-1}_{p}z^q_p + z^q_p z^{q+1}_{p} ) \Big\} \right)    \nonumber \\
            &\times  \exp \bigg\{ \tfrac{1}{T} \Big{(}   \displaystyle\sum_{j=1}^{q-1} \displaystyle\sum_{i=1}^{n-1}  x^j_i x^j_{i+1}  +  \displaystyle\sum_{i=1}^{p-2}  x^q_i x^q_{i+1} +  x^q_{p-1} z^q_{p}   - z^q_{p} y^q_{p+1} +  \displaystyle\sum_{i=p+1}^{n-1}  y^q_i y^q_{i+1}          \nonumber \\
            & +\displaystyle\sum_{j=q+1}^{n} \displaystyle\sum_{i=1}^{n-1}  y^j_i y^j_{i+1}  +   \displaystyle\sum_{i=1}^{p-1} \big{(} \displaystyle\sum_{j=1}^{q-2}  x^j_i x^{j+1}_{i}  +    x^{q-1}_i x^{q}_{i} + x_i^q z_i^{q+1} - z^{q}_i z^{q+1}_{i} +  z^{q}_i y^{q+1}_{i}   \nonumber \\
 &+   \displaystyle\sum_{j=q+1}^{n-1}  y^j_i y^{j+1}_{i}     \big{)}  +  \displaystyle\sum_{j=1}^{q-2}  x^j_p x^{j+1}_{p} +    x^{q-1}_p z^{q}_{p} - z^{q-1}_p z^{q}_{p} -  z^{q-1}_p z^{q}_{p} - z^{q}_p y^{q+1}_{p}   +   \displaystyle\sum_{j=q+1}^{n-1}  y^j_p y^{j+1}_{p}   \nonumber \\        
            &+  \displaystyle\sum_{i=p+1}^{n} \big{(}  \displaystyle\sum_{j=1}^{q-2}  x^j_i x^{j+1}_{i}  +   x^{q-1}_i z^{q}_{i} - z^{q-1}_i z^{q}_{i}  +  z^{q-1}_i y^{q}_{i} +  y^q_i y^{q+1}_{i} +   \displaystyle\sum_{j=q+1}^{n-1}  y^j_i y^{j+1}_{i}    \big{)}  \Big{)} \bigg\}. \label{config-1} 
    \end{align}
Following the line of arguments of Shiu and Chen (see \cite{Shiu chen}) we define the configurations
$$x\oplus y=\begin{pmatrix}
   x^1_1 & \cdots   & \cdots  & x^1_p& \cdots &\cdots&x^1_n  \\
      \vdots  &     &  &\vdots &  & &\vdots   \\ 
   x^{q-1}_1 & \cdots  & x^{q-1}_{p-1}& x^{q-1}_p & x^{q-1}_{p+1} &\cdots&   x^{q-1}_n \\    
    x^q_1 &\cdots  & x^q_{p-1}& z^q_p & y^q_{p+1} &\cdots&   y^{q}_n\\   
      y^{q+1}_1 & \cdots  & y^{q+1}_{p-1}& y^{q+1}_p & y^{q+1}_{p+1} &\cdots&   y^{q+1}_n \\   
  \vdots    &   & & \vdots &   && \vdots   \\ 
   y^n_1& \cdots &  \cdots & y^n_p & \cdots  &\cdots& y^n_n
\end{pmatrix} $$
 and  $x\ominus y$ with the same expression expect that in the position $(p, q)$ we have $-z^p_q$. Equation (\ref{config-1})  becomes:
\begin{align} 
Q(e)^{-1}  \pi(x) \pi(y) &= n^2 \pi(x\oplus y)  \exp\bigg\{  \tfrac{1}{T}  \Big{(}  -2( z^q_{p-1} z^q_{p} + z^{q-1}_p z^{q}_{p}   ) + \displaystyle\sum_{i=1}^{p-1}  - z^{q}_i z^{q+1}_{i}  +  \displaystyle\sum_{i=p+1}^{n}   - z^{q-1}_i z^{q}_{i}   \nonumber  \\
  &+ \displaystyle\sum_{i=1}^{p-1}  \big{(}   x_i^q z_i^{q+1} +  z^{q}_i y^{q+1}_{i}  - x^q_i y^{q+1}_{i} \big{)} -2(  z^q_{p} y^q_{p+1} + z^q_p y^{q+1}_{p} )   +\displaystyle\sum_{i=p+1}^{n}  (    x^{q-1}_i z^{q}_{i} +  z^{q-1}_i y^{q}_{i}  \nonumber \\
  & -  x^{q-1}_i y^{q}_{i}    \big{)}    \Big{)} \bigg\} + n^2 \pi(x\ominus y)  \exp \bigg\{ \tfrac{1}{T}  \Big{(}   2(z^q_p z^q_{p+1} + z^q_p z^{q+1}_{p}   )  +  \displaystyle\sum_{i=1}^{p-1}    - z^{q}_i z^{q+1}_{i}         \nonumber \\
  & -  \displaystyle\sum_{i=p+1}^{n}   z^{q-1}_i z^{q}_{i}  +  \displaystyle\sum_{i=1}^{p-1}  \big{(}   x_i^q z_i^{q+1} +  z^{q}_i y^{q+1}_{i}  - x^q_i y^{q+1}_{i} \big{)} +2(  x^q_{p-1} z^{q}_{p} + x^{q-1}_p z^{q}_{p}   )    \nonumber \\ 
  &+ \displaystyle\sum_{i=p+1}^{n}  (    x^{q-1}_i z^{q}_{i}  +  z^{q-1}_i y^{q}_{i}  - x^{q-1}_i y^{q}_{i} )  \big{)}            \Big{)} \bigg\}.  \label{config-2}
 \end{align}
 \begin{rem}
We notice that the right side in  equation (\ref{config-2}) reaches its maximum value if 
 $$-z_i^{q-1} = z_i^q = - z_i^{q+1} \,\,\,\,\,\,\,\, \mbox{for} \,\, i=1,\, \cdots,\, n .$$
  \end{rem}
Without loss of generality we can consider the situation where $ z_p^q =1$.\\
This yields:
 \begin{align} 
 Q(e)^{-1} \pi(x) \pi(y)   & \leq n^2 e^{\frac{n-1}{T}} \pi(x\oplus y)  \exp \bigg\{ \tfrac{1}{T}  \Big{(}  \displaystyle\sum_{i=1}^{p-1}  \big{(}   -x_i^q +   y^{q+1}_{i}  - x^q_i y^{q+1}_{i} \big{)}  -2(  y^q_{p+1} +      y^{q+1}_{p} )                  \nonumber  \\
  &  +\displaystyle\sum_{i=p+1}^{n}  (    x^{q-1}_i  -  y^{q}_{i} -  x^{q-1}_i y^{q}_{i}    \big{)}    \Big{)}   \bigg\}                + n^2 e^{\frac{n-1}{T}} \pi(x\ominus y)  \exp \bigg\{ \tfrac{1}{T}  \Big{(}  \displaystyle\sum_{i=1}^{p-1}  \big{(}   -x_i^q  +  y^{q+1}_{i}   \nonumber \\
&- x^q_i y^{q+1}_{i} \big{)} + 2(  x^q_{p-1}  + x^{q-1}_p   ) + \displaystyle\sum_{i=p+1}^{n}  (    x^{q-1}_i -  y^{q}_{i}  - x^{q-1}_i y^{q}_{i} )  \big{)} \Big{)}\bigg\} .\nonumber 
 \end{align}  
 Then we obtain:
 \begin{align} 
Q(e)^{-1} \sum_{\gamma_{xy} \ni e} |\gamma_{xy}| \pi(x) \pi(y)  & \leq n^4 e^{\frac{n-1}{T}} \sum_{\gamma_{xy} \ni e} \pi(x\oplus y)  \exp \bigg\{ \tfrac{1}{T}  \Big{(}  \displaystyle\sum_{i=1}^{p-1}  \big{(}   -x_i^q +   y^{q+1}_{i}  - x^q_i y^{q+1}_{i} \big{)}    \nonumber  \\
&-2(  y^q_{p+1} +      y^{q+1}_{p} ) +\displaystyle\sum_{i=p+1}^{n}  (    x^{q-1}_i  -  y^{q}_{i} -  x^{q-1}_i y^{q}_{i}    \big{)}    \Big{)} \bigg\}     \nonumber  \\
  &    + n^4 e^{\frac{n-1}{T}} \sum_{\gamma_{xy} \ni e} \pi(x\ominus y)  \exp \bigg\{  \tfrac{1}{T}  \Big{(}\displaystyle\sum_{i=1}^{p-1}  \big{(}   -x_i^q  +  y^{q+1}_{i}  - x^q_i y^{q+1}_{i} \big{)} \nonumber \\
&+2(  x^q_{p-1}  + x^{q-1}_p   ) + \displaystyle\sum_{i=p+1}^{n}  \big{(}    x^{q-1}_i -  y^{q}_{i}  - x^{q-1}_i y^{q}_{i} \big{)}   \Big{)} \bigg\} .\nonumber 
 \end{align}  
Like in Shiu and Chen  (see \cite{Shiu chen}) we notice that
 $$\bigcup_{(x,y) \,: \,\, \gamma_{xy} \ni e} \{ x\oplus y, \,\, x\ominus y  \} = \chi.$$
With new notation $w$ for the elements from $ \chi$ the previous expression becomes:
 \begin{align} 
Q(e)^{-1} \sum_{\gamma_{xy} \ni e} |\gamma_{xy}| \pi(x) \pi(y)  & \leq n^4 e^{\frac{n-1}{T}} \sum_{w \in \chi: w^p_q=1} \pi(w)  \exp \bigg\{ \tfrac{1}{T}  \Big{(} \displaystyle\sum_{i=1}^{p-1}  \big{(}   -w_i^q +   w^{q+1}_{i}  - w^q_i w^{q+1}_{i} \big{)}    \nonumber  \\
&- 2(  w^q_{p+1} +      w^{q+1}_{p} ) +\displaystyle \sum_{i=p+1}^{n}  (    w^{q-1}_i  -  w^{q}_{i} -  w^{q-1}_i w^{q}_{i}    \big{)}    \Big{)}    \bigg\}  \nonumber  \\
 &+  n^4 e^{\frac{n-1}{T}} \sum_{w \in \chi: \, w^q_p=-1} \pi(w)  \exp \bigg\{  \tfrac{1}{T}  \Big{(}  \displaystyle\sum_{i=1}^{p-1}  \big{(}   -w_i^q  +  w^{q+1}_{i}  - w^q_i w^{q+1}_{i} \big{)} \nonumber \\
&+ 2(  w^q_{p-1}  + w^{q-1}_p   ) + \displaystyle \sum_{i=p+1}^{n}  (    w^{q-1}_i -  w^{q}_{i}  - w^{q-1}_i w^{q}_{i} )  \big{)} \Big{)} \bigg\}. \nonumber 
\end{align} 
If we look to the nearest neighbors of the site $(p, q)$, for each configuration $w_{+} \in \{w \in \chi:\,\, w_p^q=1 \}$ there exist a unique configuration $w_{-} \in \{w\in \chi:\,\, w_p^q=-1 \}$ such that:
 $$ \pi(w_{+}) = \pi(w_{-} )\exp \Big\{ \tfrac{2}{T} (w_p^{q-1} + w_{p}^{q+1}+ w_{p-1}^{q} + w_{p+1}^{q}) \Big\}.$$
This yields the following equality:
\begin{align}
\displaystyle \sum_{w \in \chi: \, w^p_q=1} \pi(w) \exp \bigg\{ \tfrac{-2}{T}  \big{(}   w^q_{p+1}+ w^{q+1}_{p} \big{)}\bigg\} &=   \displaystyle \sum_{w \in \chi: \, w^p_q=-1} \pi(w) \exp \bigg\{ \tfrac{2}{T}  \big{(} w^q_{p-1}  +  w^{q-1}_p \big{)} \bigg\}. \label{symetry-1}  
\end{align}
Then, it follows that: 
 \begin{align} 
  Q(e)^{-1} \sum_{\gamma_{xy} \ni e} |\gamma_{xy}| \pi(x) \pi(y) &\leq 2 n^4 e^{\frac{1}{T}(n-1)}  \sum_{w \in \chi: w^p_q=1} \pi(w)  \exp \bigg\{ \tfrac{1}{T}  \Big{(}  \displaystyle\sum_{i=1}^{p-1}  \big{(}   -w_i^q +   w^{q+1}_{i}  - w^q_i w^{q+1}_{i} \big{)}    \nonumber  \\
&- 2(  w^q_{p+1} +      w^{q+1}_{p} ) +\displaystyle \sum_{i=p+1}^{n}  (    w^{q-1}_i  -  w^{q}_{i} -  w^{q-1}_i w^{q}_{i}    \big{)}    \Big{)} \bigg\}.     \label{Eq_somme}
  \end{align}
 \begin{rem}
 We obtain the same result by considering the case where $z_p^q=-1$.
 \end{rem} 
We turn now to  remaining cases where the site $(p, q)$ lies on the boundary of the square:\\
i) For $(p, q)=(1, 1)$ we have that  $x$ and $e^-$ coincide and $$P( e^-, e^+) =\frac{1}{n^2(1 + e^{\frac{2}{T}(z^1_{1}z^1_2 + z^1_1 z^2_{1})})}.$$
 This yields:
\begin{align} 
  Q(e)^{-1} \pi(x) \pi(y) &= n^2(1 + e^{\frac{2}{T}(z^1_{1}z^1_2 + z^1_1 z^2_{1})})  \pi(y)  \nonumber \\
         &\leq  n^2(1 + e^{\frac{4}{T}})  \pi(y) \nonumber
 \end{align}
and it follows:
 \begin{align} 
  Q(e)^{-1} \sum_{\gamma_{xy} \ni e} |\gamma_{xy}| \pi(x) \pi(y) &\leq n^4(1 + e^{\frac{4}{T}}) \sum_{y \in \chi :\, y_1^1=1 }  \pi(y)\nonumber \\  
  &= \frac{n^4}{2}(1 + e^{\frac{4}{T}}). \label{Resultat2}
   \end{align}
For  $(p, q)=(n,n)$ we notice that $y$ and $e^+$ coincide and we obtain the same kind of  result as  for $(p, q)=(1, 1)$.\\
ii) For $(p, q)=(1, n)$, we have from lemma $1.1$ $$P( e^-, e^+) = \frac{1}{n^2(1 + e^{\frac{2}{T}(  z^{1}_{n-1}z^{1}_{n} + z^{1}_{n}z^{2}_{n} )})}.$$ According to the selection of paths in section $1.1$ for any pair $(x, y)$ where  $ \gamma_{xy} \ni e$ the configurations $x$ and $y$  must have the following form:
$$\begin{array}{cc}
x=\begin{pmatrix}
   x^1_1 & \cdots   & x^{1}_{n-1}&z^1_n  \\
    z^{2}_1 & \cdots  &  \cdots&   z^{2}_n \\
      \vdots  &     &  &\vdots   \\ 
   \vdots    &   & &   \vdots   \\ 
   z^n_1& \cdots  &\cdots& z^n_n
\end{pmatrix} 
\,\,\,\,\,\,\,\,& \mbox{and} \,\,\,\,\,\,\,\,
y=\begin{pmatrix}
   z^1_1 & \cdots   & z^{1}_{n-1}&-z^1_n  \\
    y^{2}_1 & \cdots  &  \cdots&   y^{2}_n \\
      \vdots  &     &  &\vdots   \\ 
   \vdots    &   & &   \vdots   \\ 
   y^n_1& \cdots  &\cdots& y^n_n
\end{pmatrix} \end{array}.$$
As in the case where $(p, q)= \{2, \cdots, n-1\}^2$ we find: 
\begin{align}
 Q(e)^{-1} \pi(x) \pi(y) &= \frac{ n^2}{Z_T}(1 + e^{\frac{2}{T}(  z^{1}_{n-1}z^{1}_{n} + z^{1}_{n}z^{2}_{n} )})  \exp  \bigg\{ \tfrac{1}{T} \Big{(} \displaystyle\sum_{i=1}^{n-2}   x^{1}_i x^{1}_{i+1}+ x^{1}_{n-1} z^{1}_{n} - 2 z^{1}_{n-1} z^{1}_{n}       \nonumber \\
 &+\displaystyle\sum_{j=2}^{n} \displaystyle\sum_{i=1}^{n-1}  y^j_i y^j_{i+1} +   \displaystyle\sum_{i=1}^{n-1} \big{(} x^1_i z^{2}_{i} + z^1_i y^{2}_{i}  +   \displaystyle\sum_{j=2}^{n-1}  y^j_i y^{j+1}_{i}   - z^1_i z^{2}_{i}  \big{)} - z^1_n y^{2}_{n} + \displaystyle\sum_{j=2}^{n}   y^j_n y^{j+1}_{n}  \Big{)} \Bigg\}. \label{config-3}
 \end{align}
Following the line of arguments from Shiu and Chen (see \cite{Shiu chen}) we define configurations:
$$\begin{array}{cc}
x\oplus y=\begin{pmatrix}
   x^1_1 & \cdots   & x^{1}_{n-1}&z^1_n  \\
    y^{2}_1 & \cdots  &  \cdots&   y^{2}_n \\
      \vdots  &     &  &\vdots   \\ 
   \vdots    &   & &   \vdots   \\ 
   y^n_1& \cdots  &\cdots& y^n_n
\end{pmatrix} 
\,\,\,\,\,\,\,\,& \mbox{and} \,\,\,\,\,\,\,\,
x \ominus y=\begin{pmatrix}
   x^1_1 & \cdots   & x^{1}_{n-1}&-z^1_n  \\
    y^{2}_1 & \cdots  &  \cdots&   y^{2}_n \\
      \vdots  &     &  &\vdots   \\ 
   \vdots    &   & &   \vdots   \\ 
   y^n_1& \cdots  &\cdots& y^n_n
\end{pmatrix} \end{array}.$$
Then equation (\ref{config-3}) becomes: 
\begin{align}
 Q(e)^{-1} \pi(x) \pi(y) &= n^2 \pi(x\oplus y)  \exp \bigg\{ \tfrac{1}{T} \Big{(} -2 (z_{n-1}^1 z_n^1 + z_{n}^{1} y_{n}^{2}  ) + \displaystyle \sum_{i=1}^{n-1} \big{(}  x_{i}^1 z_i^2    +    z_{i}^{1} y_{i}^{2}   -  z_{i}^{1} z_{i}^{2} -  x_{i}^{1} y_{i}^{2}    \big{)}    \Big{)} \bigg\} \nonumber \\
 &+ n^2 \pi(x\ominus y)  \exp \bigg\{ \tfrac{1}{T} \Big{(} 2(z_n^1 z_n^2 + x_{n-1}^1z_n^1) + \displaystyle \sum_{i=1}^{n-1}   \big{(}  x_{i}^1 z_i^2    +    z_{i}^{1} y_{i}^{2}   -  z_{i}^{1} z_{i}^{2} -  x_{i}^{1} y_{i}^{2}    \big{)}     \Big{)} \bigg\}.   \label{config-4}
\end{align}
\begin{rem}
We notice that equation (\ref{config-4}) reaches its  maximum if the following two points hold:
 \begin{itemize}
\item $z_i^1 = - z_i^2=-1$  for  $i \leq n-1$.   
\item $z_n^1 = z_n^2 =+1$ and $z_{n-1}^1 =- z_n^1 = -1 $.
\end{itemize}
\end{rem}
This yields:
\begin{align}
 Q(e)^{-1} \pi(x) \pi(y) &\leq n^2 \pi(x\oplus y)  \exp \bigg\{ \tfrac{1}{T} \Big{(} 2 (1 - y_{n}^{2}  ) + \displaystyle \sum_{i=1}^{n-1} \big{(}  x_{i}^1  - y_{i}^{2}  +1 -  x_{i}^{1} y_{i}^{2}    \big{)}    \Big{)} \bigg\} \nonumber \\
 &+ n^2 \pi(x\ominus y)  \exp \bigg\{ \tfrac{1}{T} \Big{(} 2(1+ x_{n-1}^1) + \displaystyle \sum_{i=1}^{n-1}   \big{(}  x_{i}^1     - y_{i}^{2}   +1 -  x_{i}^{1} y_{i}^{2}    \big{)}     \Big{)} \bigg\}.  \label{config-5}
\end{align}
Summing over  all pairs $(x, y)$ where the path $\gamma_{xy}$ passes through the edge $e$ of equation (\ref{config-5}) gives:
\begin{align}
 Q(e)^{-1} \sum_{\gamma_{xy} \ni e} |\gamma_{xy}| \pi(x) \pi(y) &\leq n^4 e^{\frac{n-1}{T}} \sum_{\gamma_{xy} \ni e} \pi(x\oplus y)  \exp\bigg\{ \tfrac{1}{T} \Big{(} 2 (1 - y_{n}^{2}  ) + \displaystyle \sum_{i=1}^{n-1} \big{(}  x_{i}^1  - y_{i}^{2}  -  x_{i}^{1} y_{i}^{2}    \big{)}    \Big{)} \bigg\} \nonumber \\
 &+ n^4 e^{\frac{n-1}{T}} \sum_{\gamma_{xy} \ni e} \pi(x\ominus y)  \exp \bigg\{ \tfrac{1}{T} \Big{(} 2(1+ x_{n-1}^1) + \displaystyle \sum_{i=1}^{n-1}   \big{(}  x_{i}^1     - y_{i}^{2}  -  x_{i}^{1} y_{i}^{2}    \big{)}     \Big{)}   \bigg\}\nonumber \\
 &= n^4 e^{\frac{n-1}{T}} \sum_{w \in \chi:\,\, w_n^1=1} \pi(w)  \exp \bigg\{ \tfrac{1}{T} \Big{(} 2 (1 - w_{n}^{2}  ) + \displaystyle \sum_{i=1}^{n-1} \big{(}  w_{i}^1  - w_{i}^{2}  -  w_{i}^{1} w_{i}^{2}    \big{)}    \Big{)} \bigg\} \nonumber \\
 &+ n^4 e^{\frac{n-1}{T}} \sum_{w \in \chi:\,\, w_n^1=-1} \pi(w)  \exp \bigg\{ \tfrac{1}{T} \Big{(} 2(1+ w_{n-1}^1) + \displaystyle \sum_{i=1}^{n-1}   \big{(}  w_{i}^1     - w_{i}^{2}   -  w_{i}^{1} w_{i}^{2}    \big{)}     \Big{)}   \bigg\}\nonumber \\
 &= 2 n^4 e^{\frac{n-1}{T}} \sum_{w \in \chi:\,\, w_n^1=1} \pi(w)  \exp \bigg\{ \tfrac{1}{T} \Big{(} 2 (1 - w_{n}^{2}  ) + \displaystyle \sum_{i=1}^{n-1} \big{(}  w_{i}^1  - w_{i}^{2} -  w_{i}^{1} w_{i}^{2}    \big{)}    \Big{)} \bigg\}. \nonumber
 \end{align}
 The last equality is given by a symmetry argument similar to the one introduced in equation (\ref{symetry-1}).
 For  $(p, q)=(n,1)$ we use the same technique and we obtain a similar result.\\
iii) For $(p ,q) \in \{2 \cdots n-1\} \times \{1\}$, we have from lemma $1.1$ $$P( e^-, e^+)=\frac{1}{n^2(1 + e^{\frac{2}{T}(  z^{q-1}_{1}z^{q}_{1} + z^{q}_{1}z^{q+1}_{1} + z^{q}_{1}z^{q}_{2})})}.$$ 
As before, a pair of configuration $(x, y)$ where $\gamma_{xy} \ni e$ must have the following form:
$$\begin{array}{cc}
x=\begin{pmatrix}
   x^1_1 & \cdots   & \cdots&x^1_n  \\
      \vdots  &     &  &\vdots   \\ 
   x^{q-1}_1 & \cdots  &  \cdots&   x^{q-1}_n \\    
    z^q_1 &\cdots  & \cdots&   z^{q}_n\\   
      z^{q+1}_1 & \cdots  &\cdots&   z^{q+1}_n \\   
  \vdots    &   & &   \vdots   \\ 
   z^n_1& \cdots  &\cdots& z^n_n
\end{pmatrix} 
\,\,\,\,\,\,\,\,& \mbox{and} \,\,\,\,\,\,\,\,
y=\begin{pmatrix}
   z^1_1 & \cdots   & \cdots&z^1_n  \\
      \vdots  &     &  &\vdots   \\ 
   z^{q-1}_1 & \cdots  &  \cdots&   z^{q-1}_n \\    
    -z^q_1 &y^q_2  & \cdots&   y^{q}_n\\   
      y^{q+1}_1 & \cdots  &\cdots&   y^{q+1}_n \\   
  \vdots    &   & &   \vdots   \\ 
   y^n_1& \cdots  &\cdots& y^n_n
\end{pmatrix} 
\end{array}.$$
It follows that:
\begin{align}
 Q(e)^{-1} \pi(x) \pi(y) &= \frac{ n^2}{Z_T}(1 + e^{\frac{2}{T}(  z^{q-1}_{1}z^{q}_{1} + z^{q}_{1}z^{q+1}_{1} + z^{q}_{1}z^{q}_{2})}) \nonumber \\
 &\times \exp\bigg\{ \tfrac{1}{T} \Big{(} \displaystyle\sum_{j=1}^{q-1} \displaystyle\sum_{i=1}^{n-1}  x^{j}_i x^{j}_{i+1}  -  z^q_1 y^q_{2} +  \displaystyle\sum_{i=2}^{n-1}  y^q_i y^q_{i+1}   +\displaystyle\sum_{j=q+1}^{n} \displaystyle\sum_{i=1}^{n-1}  y^j_i y^j_{i+1}  \nonumber \\ 
&+  \displaystyle\sum_{j=1}^{q-2}  x^j_1 x^{j+1}_{1} +    x^{q-1}_1 z^{q}_{1}  -    z^{q-1}_1 z^{q}_{1} - z^q_1 y^{q+1}_{1}   +   \displaystyle\sum_{j=q+1}^{n-1}  y^j_1 y^{j+1}_{1}  - z^{q-1}_1 z^{q}_{1}   \nonumber \\ 
&+  \displaystyle\sum_{i=2}^{n} \big{(}  \displaystyle\sum_{j=1}^{q-2}  x^j_i x^{j+1}_{i}  +   x^{q-1}_i z^{q}_{i}    +   z^{q-1}_i y^{q}_{i} +  y^q_i y^{q+1}_{i} +   \displaystyle\sum_{j=q+1}^{n-1}  y^j_i y^{j+1}_{i}      -   z^{q-1}_i z^{q}_{i}  \big{)}   \Big{)} \bigg\}. \label{config-6} 
\end{align}
Following the line of arguments of Shiu and Chen  (see \cite{Shiu chen}) we define the configurations:
$$\begin{array}{cc}
x \oplus y=\begin{pmatrix}
   x^1_1 & \cdots   & \cdots&x^1_n  \\
      \vdots  &     &  &\vdots   \\ 
   x^{q-1}_1 & \cdots  &  \cdots&   x^{q-1}_n \\    
    z^q_1 & y^q_2 & \cdots&   y^{q}_n\\   
      y^{q+1}_1 & \cdots  &\cdots&   y^{q+1}_n \\   
  \vdots    &   & &   \vdots   \\ 
   y^n_1& \cdots  &\cdots& y^n_n
\end{pmatrix}
\,\,\,\,\,\,\,\,& \mbox{and}\,\,\,\,\,\,\,\,
x\ominus y=\begin{pmatrix}
   x^1_1 & \cdots   & \cdots&x^1_n  \\
      \vdots  &     &  &\vdots   \\ 
   x^{q-1}_1 & \cdots  &  \cdots&   x^{q-1}_n \\    
    -z^q_1 &y^q_2  & \cdots&   y^{q}_n\\   
      y^{q+1}_1 & \cdots  &\cdots&   y^{q+1}_n \\   
  \vdots    &   & &   \vdots   \\ 
   y^n_1& \cdots  &\cdots& y^n_n
\end{pmatrix} 
\end{array}.$$
Equation (\ref{config-6}) becomes:
\begin{align}
 \frac{\pi(x) \pi(y)}{Q(e)} &= n^2 \pi(x \oplus y)   \exp\bigg\{ \tfrac{1}{T} \big{(}  - 2 (z^{q}_1 y^{q}_{2} + z^{q-1}_1 z^{q}_{1} + z_1^qy_1^{q+1}) + \displaystyle\sum_{i=2}^{n} \big{(}    x^{q-1}_i z^{q}_{i}    +   z^{q-1}_i y^{q}_{i}    \nonumber \\
&  -   z^{q-1}_i z^{q}_{i} -   x^{q-1}_i y^{q}_{i}  \big{)}   \big{)} \bigg\}    +   n^2 \pi(x \ominus y)   \exp\bigg\{ \tfrac{1}{T} \big{(} 2(  z^{q}_{1}z^{q+1}_{1} +  z^{q}_{1}z^{q}_{2} +    x^{q-1}_1 z^{q}_{1} )  \nonumber \\
&+  \displaystyle\sum_{i=2}^{n} \big{(}   x^{q-1}_i z^{q}_{i}    +   z^{q-1}_i y^{q}_{i}      -   z^{q-1}_i z^{q}_{i} -   x^{q-1}_i y^{q}_{i}  \big{)}   \big{)} \bigg\}.  \label{config-7}  
\end{align}
\begin{rem}
The maximum in equation (\ref{config-7})  is reached once the following conditions hold:
\begin{itemize}
\item $-z_i^{q-1} = z_i^q = z_i^{q+1}$ \,\,\,\,\,\, for $i=2, \cdots, n$.
\item $-z_1^{q-1} = z_1^q = z_1^{q+1} =1 $ and  $ z_1^{q} = z_2^q =1$.
\end{itemize}
\end{rem}
Then equation (\ref{config-7}) becomes:
\begin{align}
 Q(e)^{-1} \pi(x) \pi(y) &\leq n^2 \pi(x \oplus y)   \exp \bigg\{ \tfrac{1}{T} \Big{(}  - 2 (y^{q}_{2} -1 + y_1^{q+1})   + \displaystyle\sum_{i=2}^{n} \big{(}    x^{q-1}_i     -  y^{q}_{i}   +1  -   x^{q-1}_i y^{q}_{i}  \big{)}   \Big{)}  \Bigg\}\nonumber \\
&+   n^2 \pi(x \ominus y)   \exp \bigg\{ \tfrac{1}{T} \Big{(} 2(  1 +  1  +    x^{q-1}_1  ) +  \displaystyle\sum_{i=2}^{n} \big{(}   x^{q-1}_i    - y^{q}_{i}  +1 -   x^{q-1}_i y^{q}_{i}  \big{)}   \Big{)} \bigg\}.     \nonumber 
\end{align}
As in Shiu and Chen (see \cite{Shiu chen}) we notice that $$ \bigcup_{(x, y):\, \gamma_{xy} \ni e} \Big\{ x\oplus y, \, x\ominus y \Big\} = \chi. $$
 In the previous expression we change the notation of an element   $x\oplus y$ or $x\ominus y$ to $w$ with $w_p^1=+1$ or $w_p^1=-1$ respectively. This yields:
\begin{align}
 Q(e)^{-1} \sum_{\gamma_{xy} \ni e} |\gamma_{xy}| \pi(x) \pi(y) &\leq n^4 \sum_{\gamma_{xy} \ni e}   \pi(x \oplus y)   \exp \bigg\{ \tfrac{1}{T} \Big{(}   2 (1-y^{q}_{2}  - y_1^{q+1})  + \displaystyle\sum_{i=2}^{n} \big{(}    x^{q-1}_i     -  y^{q}_{i}         \nonumber \\
&+1 -   x^{q-1}_i y^{q}_{i}  \big{)}   \Big{)} \bigg\}   +   n^4 \sum_{\gamma_{xy} \ni e} \pi(x \ominus y)   \exp \bigg\{ \tfrac{1}{T} \Big{(} 2(  1 +  1  +    x^{q-1}_1  )   \nonumber \\
&+  \displaystyle\sum_{i=2}^{n} \big{(}   x^{q-1}_i   - y^{q}_{i}  +1  -   x^{q-1}_i y^{q}_{i}  \big{)}   \Big{)}  \bigg\}        \nonumber \\
&= n^4 e^{\frac{n+1}{T}} \sum_{w\in \chi: \,w_p^1=1 }  \pi(w)   \exp\bigg\{ \tfrac{1}{T} \Big{(} - 2 (w^{q}_{2} +w_1^{q+1})   + \displaystyle\sum_{i=2}^{n} \big{(}    w^{q-1}_i      \nonumber \\
&   -  w^{q}_{i} -   w^{q-1}_i w^{q}_{i}  \big{)}  \Big{)}  \bigg\}    +  n^4 e^{\frac{n+1}{T}}  \sum_{w\in \chi: \,w_p^1=-1}  \pi(w)   \exp \bigg\{ \tfrac{1}{T} \Big{(} 2(  1  +    w^{q-1}_1  )  \nonumber \\
&+  \displaystyle\sum_{i=2}^{n} \big{(}   w^{q-1}_i    - w^{q}_{i}  -   w^{q-1}_i w^{q}_{i}  \big{)}   \Big{)}   \bigg\} . \nonumber 
\end{align}
We obtain the same result for  $(p, q)\in\{2\cdots n-1\}\times \{n\}$. With similar techniques we can do the same calculus for $(p, q)\in \{1\} \times \{ 2 \cdots n-1\}$ then we obtain:
\begin{align}
 Q(e)^{-1} \sum_{\gamma_{xy} \ni e} |\gamma_{xy}| \pi(x) \pi(y) &\leq n^4 e^{\frac{n+1}{T}} \sum_{w\in \chi: \,w_p^1=1 }  \pi(w)   \exp\bigg\{ \tfrac{1}{T} \Big{(} - 2 (w^{1}_{p+1} +w_p^{2})   + \displaystyle\sum_{i=1}^{n-1} \big{(}    w^{1}_i     \nonumber \\
&  -  w^{2}_{i} -   w^{1}_i w^{2}_{i}  \big{)}  \Big{)}  \bigg\}      +  n^4 e^{\frac{n+1}{T}}  \sum_{w\in \chi: \,w_p^1=-1}  \pi(w)   \exp \bigg\{ \tfrac{1}{T} \Big{(} 2(  1+    w^{1}_{p-1})  \nonumber \\
&  +  \displaystyle\sum_{i=1}^{n-1} \big{(}   w^{1}_i    - w^{2}_{i}  -   w^{1}_i w^{2}_{i}  \big{)}   \Big{)}   \bigg\} . \nonumber
\end{align} 
A similar result is obtained in the case where $(p, q)\in \{n\} \times \{ 2 \cdots n-1\}$.

We turn now to give the proof of the main theorem of this paper:
\subsection{Proof of  theorem $2$}  
  Since for all sites $(p, q) \in \{ 1, \cdots, n \}^2$ we have $w_p^q = \pm 1$ it follows that
$$\displaystyle \max_{i\in \{1, \cdots, n\}} \{ w^{q-1}_i  -  w^{q}_{i} -  w^{q-1}_i w^{q}_{i} \} = \displaystyle \max_{i\in \{1, \cdots, n\}} \{   -w_i^q +   w^{q+1}_{i}  - w^q_i w^{q+1}_{i} \} =3.$$
Analyzing the worst cases yield the following inequalities \\
 i) $2 (1 - w_{n}^{2}  )+ \displaystyle \sum_{i=1}^{n-1} \big{(}  w_{i}^1  - w_{i}^{2}  -  w_{i}^{1} w_{i}^{2}    \big{)} \leq 3n+1.$ \\
ii) $-2( w_2^q + w_1^{q+1}) + \displaystyle\sum_{i=2}^{n} \big{(}    w^{q-1}_i  - w^{q}_{i}  -   w^{q-1}_i w^{q}_{i}  \big{)} \leq 3n+1.$\\
iii) $2( 1 + w_1^{q-1}) + \displaystyle\sum_{i=2}^{n} \big{(}    w^{q-1}_i  - w^{q}_{i}  -   w^{q-1}_i w^{q}_{i}  \big{)} \leq 3n+1.$\\
iv) $ \displaystyle\sum_{i=1}^{p-1}  \big{(}   -w_i^q +   w^{q+1}_{i}  - w^q_i w^{q+1}_{i} \big{)} - 2(  w^q_{p+1} +w^{q+1}_{p} ) +\displaystyle \sum_{i=p+1}^{n}  ( w^{q-1}_i  -  w^{q}_{i} -  w^{q-1}_i w^{q}_{i}    \big{)} \leq 3n+1. $\\ \\
By the fact that $$\displaystyle\sum_{w \in \chi:\,\, w_p^q =+1} \pi(w) = \displaystyle\sum_{w \in \chi:\,\, w_q^q =-1} \pi(w) = \frac{1}{2}$$ we conclude that $$ \kappa \leq n^4 \exp \Big\{ \tfrac{2}{T} (2n+1) \Big\} $$ which finishes the proof.
\section*{Acknowledgements}  
The authors would  like to thank Professor Mondher Damak from University of Sfax for his help with this paper.\\
This work was supported by "Direction G\'en\'erale de la Recherche Scientifique de la R\'epublique Tunisienne, Minist\`ere de l'Enseignement Sup\'erieur et de la Recherche Scientifique" and "Service de Coop\'eration et d'Action Culturelle de l'Ambassade de la R\'epublique Fran\c{c}aise en Tunisie", PHC Utique project Number 16G1505.

\end{document}